\documentclass[conference,a4paper]{IEEEtran}
\usepackage{diagbox,setspace}
\usepackage{xcolor}
\usepackage{balance}
\usepackage{cite}
\ifCLASSINFOpdf
  \usepackage[pdftex]{graphicx}
\else
  \usepackage[dvips]{graphicx}
\fi
\usepackage[cmex10]{amsmath}
\usepackage{mathtools,tikz}
\usepackage{pifont}
\ifCLASSOPTIONcompsoc
 \usepackage[caption=false,font=normalsize,labelfont=sf,textfont=sf]{subfig}
\else
 \usepackage[caption=false,font=footnotesize]{subfig}
\fi
\usepackage[english]{babel}
\usepackage{ifluatex}
\usepackage[utf8]{\ifluatex lua\fi inputenc}
\usepackage[T1]{fontenc}
\usepackage{textcomp}
\usepackage[scaled=0.91]{helvet}
\usepackage[cmintegrals,varvw,]{newtxmath}
\usepackage{shellesc}
\usepackage{bm}
\usepackage{siunitx,xcolor,booktabs}

\newcommand*\mat[1]{\textsf{\textit{\textbf#1}}}

\renewcommand*\vec[1]{\textsf{\textit{\textbf#1}}}

\DeclarePairedDelimiterX\braket[2]{\langle}{\rangle}{#1 , #2}

\usepackage[a4paper,top=18mm,left=14.3mm,right=14.3mm,bottom=44.3mm]{geometry}%
\begin{document}
\title{Improved Discretization of the Full First-Order Magnetic Field Integral Equation }

\author{\IEEEauthorblockN{
Jonas Kornprobst\IEEEauthorrefmark{1},
Alexander Paulus\IEEEauthorrefmark{1}, and   
Thomas F. Eibert\IEEEauthorrefmark{1}
}
\IEEEauthorblockA{\IEEEauthorrefmark{1}
Chair of High-Frequency Engineering\\Department of Electrical and Computer Engineering\\Technical University of Munich\\80290 Munich, Germany\\e-mail: \{j.kornprobst;eibert\}@tum.de}
}

\maketitle

\begin{abstract}
The inaccuracy of the classical magnetic field integral equation (MFIE) is a long-studied problem.
We investigate one of the potential approaches to solve the accuracy problem: higher-order discretization schemes.
While these are able to offer increased accuracy, we demonstrate that the accuracy problem may still be present. 
We propose an advanced scheme based on a weak-form discretization of the identity operator which is able to improve the high-frequency MFIE accuracy considerably\,---\,without any significant increase in computational effort or complexity. 
\end{abstract}

\vskip0.5\baselineskip

\section{Introduction} %
The method of moments is a common choice for solving scattering and radiation problems\cite{harrington1968field}. 
Considering perfect electrically conducting (PEC) objects, there are two basic choices for the discretization of the relevant integral operators.
On the one hand, the electric field integral equation (EFIE) is known to be quite accurate even with low-order discretization schemes.
It is applicable to open and closed bodies, but it suffers from ill-conditioning at low frequencies and for dense discretizations.
On the other hand, the magnetic field integral equation (MFIE) is applicable to closed bodies only and offers a well-conditioned system matrix regardless of the considered scenario. However, the accuracy of the MFIE is a topic of concern and of ongoing research.

Historically, the MFIE accuracy problems have been observed for lowest-order (LO) divergence-conforming discretizations, i.e., with Rao-Wilton-Glisson (RWG) basis functions used for the electric surface current densities and for the testing of the magnetic field\cite{Rao1982RWG,Rius2001solidangleMFIE,zhang2003magnetic,Ergul2004solidangle,Ergul2009Discretization}. 
These inaccuracies are in particular observed for low-frequency scenarios and for objects with sharp edges.
The problems have several (related) reasons. 
Instead of the mathematically advantageous curl-conforming testing functions in the dual space of the RWG functions, div-conforming RWGs are employed for the classical testing of the MFIE. Furthermore, a considerable part of the inaccuracies stems from the fact that the MFIE is an integral equation of the second kind. 
The involved highly-singular identity operator diminishes some of the positive averaging effects of the integrals.

The first comprehensive solution to these issues has been a mathematical approach: testing the magnetic field of RWG functions in its dual space, i.e., with Buffa-Christiansen (BC) functions\cite{chen1990electromagnetic,buffa2007dual,tong2009dual,Cools2011BCMFIE,Yan2011BCMFIE}. 
The BC-tested MFIE does not suffer from the mentioned problems anymore, however, at the cost of handling the more demanding barycentrically refined mesh.
The low-frequency breakdown of the MFIE has also been tackled by razor-blade testing\cite{bogaert2013low}. 
Since high-frequency applications are more relevant for most electromagnetic applications, curing just the high-frequency problems has in contrast attracted much more attention.
Still, there are only a few approaches which are able to provide reliable solutions.
Monopolar RWGs have been investigated some time ago\cite{ubeda2006novel}. 
However, for their combination with the EFIE it is necessary to employ volumetric zero-field testing functions due to the involved unphysical hypersingular integrals\cite{ubeda2020accurate}.
Other approaches with curl-conforming ansatz functions (e.g., with rotated RWGs) suffer from the same problem\cite{peterson2002curlMFIE,ubeda2005mfiecurlbasis,ergul2006curlMFIE}.
Going for higher-order (HO) discretizations certainly increases the MFIE accuracy and hence seems to eliminate the accuracy issues, at least at first sight~\cite{Pasi2005TLbasis,Ergul2007Linearlinear,eibert2009surface}.
Rather recently, a weak-form (WF) RWG discretization of the MFIE operator and in particular the involved identity operator was proposed, which is able to solve the high-frequency accuracy problem without introducing non-integrable sources or worsening the conditioning~\cite{kornprobst2018combined,kornprobst2018weak,kornprobst2018accurate2}.

In this paper, we show that\,---\,unlike sometimes claimed in literature\,---\,higher-order functions may not fully solve the MFIE accuracy issues. 
To do so, we consider a full first-order discretization which was considered in\cite{sun2001construction,eibert2009surface,li2014analytical} and should basically be similar to other full-first-order approaches involving Trintinalia-Ling functions~\cite{trintinalia2001first}.
In order to arrive at an accurate discretization, it may still be necessary to cure the RWG part of the identity operator discretization just as this is necessary for the LO discretization alone.

The paper is structured as follows. The MFIE and its discretization are introduced in Section II. 
Section III discusses the proposed improved discretization. 
The results in Section IV demonstrate the achieved accuracy improvements.

\section{The Magnetic Field Integral Equation} %
The equivalence principle states that the electromagnetic behavior of a PEC object may be modeled by equivalent electric and magnetic surface current densities on its surface $s$.
In order to arrive at the MFIE of PEC scattering scenarios, we choose equivalent electric  Love surface current densities on $s$ according to 
\begin{equation}
\bm j = \bm n \times \bm h
\end{equation}
with the outward unit normal vector $\bm n$ on $s$, which in turn are the basis of an equivalent scenario in free space involving the associated well-known Green's functions of the scalar Helmholtz equation $g(\bm r,\bm r')$.
Enforcing the Love condition for the radiated magnetic field on $s$ yields the MFIE
\begin{equation}
\big[\frac{1}{2}\bm{\mathcal I}+\bm{\mathcal K}\big]\bm j = \bm h^\mathrm{inc} \times \bm n
\end{equation}
with the identity operator $\bm{\mathcal I}$ and the MFIE integral operator 
\begin{equation}
\bm{\mathcal K}\{\bm j\} = - \bm n (\bm r) \times \iint_s \bm \nabla g(\bm r,\bm r') \times \bm j (\bm r) \,\mathrm{d}^2r'.
\end{equation}
The classical discretization of the MFIE with ansatz functions $\bm\beta$ (their rotated counterparts are $\bm\alpha = \bm\beta\times\bm n$) for the electric surface currents 
\begin{equation}
\bm j = \sum_{n=1}^N \bm \beta_n [\vec i]_n 
\end{equation}
then reads
\begin{equation}
\big[\frac{1}{2}\mat G_{\bm\beta,\bm\beta}+\mat K_{\bm\beta,\bm\beta}\big]\vec i = \vec h ^\mathrm{\,inc}_{\bm \beta}
\label{eq:old_mfie}
\end{equation}
with the matrix entries and the right-hand side defined in the style of inner products
\begin{equation}
[\mat G_{\bm\beta,\bm\beta}]_{mn}=  \braket{\bm\beta_m}{\bm{\mathcal I}\{\bm\beta_n \}}
\,,\quad
[\mat K_{\bm\beta,\bm\beta}]_{mn}=  \braket{\bm\beta_m}{\bm{\mathcal K}\{\bm\beta_n \}}
\,,
\end{equation}
\begin{equation}
[\vec h ^\mathrm{\,inc}_{\bm \beta}]_{m}=  \braket{\bm\beta_m}{\bm h^\mathrm{inc} \times \bm n}=\iint_s \bm\beta_m\cdot (\bm h^\mathrm{inc} \times \bm n)\,\mathrm{d}^2r
\,.
\end{equation}
The $\bm\beta$ expansion functions of the electric currents may be just $N$ RWG functions defined on all neighboring pairs of triangular mesh facets.
For a full first-order discretization, we split these functions into an RWG part with $N/2$ functions $\bm \beta_\textsc{lo}$ and the complementary extension to full first order~\cite{eibert2009surface,li2014analytical} with $N/2$ functions  $\bm \beta_\textsc{ho}$.
The following investigations are based on this hierarchical HO approach. 
For other HO ansatz functions\,---\,in particular, if they are of interpolatory nature\,\mbox{---,} it may not be that easy to identify and isolate the anisotropic influence of the RWG-alike part of the functions.

\section{Improved Discretization of the MFIE} %
The basis of the MFIE accuracy improvements is a WF discretization of the identity operator, i.e., the Gram matrix block of the RWG functions in~\eqref{eq:old_mfie} is (partially) replaced, leading to~\cite{kornprobst2018accurate2}
\begin{equation}
\big[\frac{1}{4}\mat G_{\bm\beta,\bm\beta}-\frac{1}{4}\mat G_{\bm\alpha,\bm\beta}\mat G_{\bm\beta,\bm\beta}^{\smash{-1}}\mat G_{\bm\alpha,\bm\beta}+\mat K_{\bm\beta,\bm\beta}\big]\vec i = \vec h ^\mathrm{\,inc}_{\bm \beta}
\label{eq:new_mfie}
\end{equation}
with the Gram matrix of rotated and standard RWGs
\begin{equation}
[\mat G_{\bm\alpha,\bm\beta}]_{mn}=  \braket{\bm\alpha_m}{\bm\beta_n}
\,.
\end{equation}
As seen in~\eqref{eq:new_mfie}, a weighting of 50\% of the classical RWG Gram matrix and 50\% of the WF discretization\,---\,two subsequent WF $90^\circ$ rotations and a sign inversion\,---\,is a suitable choice for pure RWG discretization.

In the implementation of~\eqref{eq:old_mfie}, it is straight-forward to split the matrices into RWG and HO parts.
In the full first-order discretization, the HO part already reduces the anisotropy of the RWG functions  to some extent.
Hence, the “correction term”  needs here a weaker weighting than 50\%.
We have chosen 20\% for the results presented in this paper.

\section{Numerical Results} %
\begin{figure}[t]
\centering
\begin{minipage}[t]{\linewidth}\centering
  \subfloat[\hspace*{1.5cm}]{%
\includegraphics{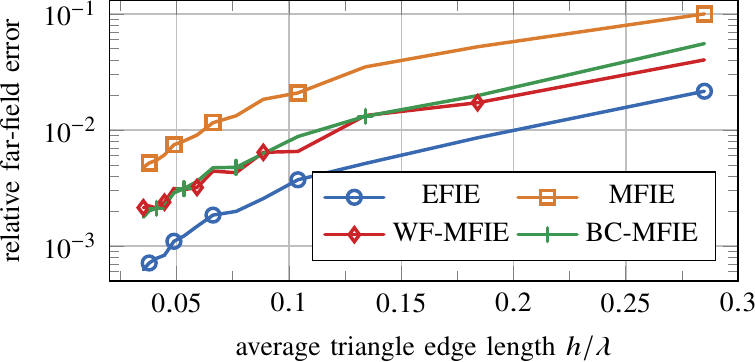}}
\end{minipage}
\\
\begin{minipage}[t]{\linewidth}\centering
  \subfloat[\hspace*{1.5cm}]{%
\includegraphics{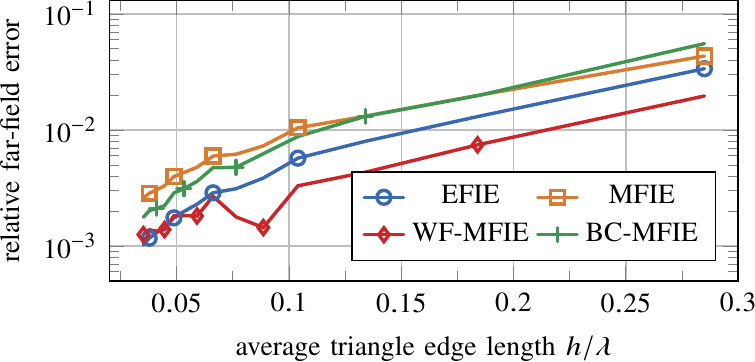}}
\end{minipage}
\caption{Mesh refinement accuracy analysis of the averaged radar-cross section of a cuboid for the comparison of classical EFIE and MFIE, as well as BC-tested MFIE and and MFIE with WF-discretized identity operator. (a) RWG LO discretization. (b) Full first-order discretization of all IEs except for the BC-MFIE from BEAST.\label{fig:mesh}}
\end{figure}

Plane-wave scattering from a square cuboid with an edge length of $0.5\lambda$ is analyzed in a mesh refinement study in Fig.~\ref{fig:mesh} for RWG and full first-order ansatz functions.
For comparison purposes, the Julia implementation of the BC-MFIE in the boundary element analysis and simulation toolkit (BEAST) was employed~\cite{kristof_cools_2020_3865280}.
The errors are calculated with respect to a $2.5$th order EFIE solution.
As it is known from literature~\cite{Pasi2005TLbasis,Ergul2007Linearlinear,eibert2009surface}, we observe that the full first-order discretization of the MFIE increases the accuracy, here to a level close to the BC-MFIE and the WF-MIFE\,---\,which are almost indistinguishable in their accuracy level.
However, the proposed HO WF-MFIE scheme is able to increase the accuracy further.
This clearly hints to the fact that the classical HO MFIE does not solve the MFIE problems but only relieves them in part.
For a fully satisfying HO MFIE, dual-space testing would certainly be required.
Nevertheless, the proposed WF scheme pushes the MFIE accuracy to EFIE levels for high-frequency scattering scenarios.
As seen in Fig.~\ref{fig:mesh}, the full first-order WF-MFIE is here even more accurate than the full-first order EFIE, mostly due to the fact that the EFIE cannot benefit from the full first-order discretization in the considered frequency range. 

\section{Conclusion} %
We have presented a WF discretization scheme for the MFIE with classical div-conforming testing functions.
The proposed scheme works with full first-order ansatz functions and is able to increase the accuracy of the MFIE considerably.
Furthermore, it becomes clear from the presented results that the standard HO MFIE does not fully solve the MFIE accuracy problems but only hides them in part.

\bibliographystyle{IEEEtran}
\bibliography{IEEEabrv,ref}

\end{document}